\documentclass{amsart}
\usepackage{almpos}

\usepackage[colorlinks=true,
citecolor=black,
linkcolor=black,
anchorcolor=black,
filecolor=black,
menucolor=black,
urlcolor=black,
pdftitle={Nilpotency, almost nonnegative curvature, and gradient flow},
pdfsubject={Differential geometry},
pdfauthor={Vitali Kapovitch, Anton Petrunin and Wilderich Tuschmann}
]{hyperref}

\begin{document}
\title[Nilpotency, almost nonnegative curvature, and gradient flow]
{Nilpotency, almost nonnegative curvature, and the gradient flow on
Alexandrov spaces}

\rm
\thanks{\it The first two authors were supported in part by NSF grants \# DMS-0204187 and \# DMS-0406482 respectively.
The last author's research was partially supported by a DFG
Heisenberg Fellowship.}
\author{Vitali Kapovitch}
\address{Vitali Kapovitch\\Department of Mathematics\\University of Toronto\\
Toronto, Ontario, M5S 2E4, Canada}\email{vtk@math.toronto.edu}
\author{Anton Petrunin }\address{Anton Petrunin\\ Department of Mathematics\\ Pennsylvania State University\\
University Park, State College, PA 16802
}\email{petrunin@math.psu.edu}
\author{Wilderich Tuschmann}\address{Wilderich Tuschmann
\\Mathematisches Seminar
\\Christian-Albrechts--Universit\"at zu Kiel\\ Ludewig-Meyn--Stra\ss{}e 4-8 \\ D-24118 Kiel,
Germany}\email{tusch@math.uni-kiel.de}

\begin{abstract}
We show that almost nonnegatively curved $m$-manifolds are, up to finite cover,
nilpotent spaces in the sense of homotopy theory
and have $C(m)$-nilpotent fundamental groups.
We also show that up to a finite cover
almost nonnegatively curved manifolds are fiber bundles
with simply connected fibers over nilmanifolds.
\end{abstract}

\maketitle

\section{Introduction}

Almost nonnegatively curved manifolds were introduced by Gromov
in the late 70s~\cite{G4},
with the most significant contributions to their study made by Fukaya and Yamaguchi in \cite{Yam} and \cite{FY}.
Building on their ideas,
in the present article
we establish several
new properties of these 
manifolds which yield, in particular,
new topological obstructions to almost nonnegative curvature.
 Our techniques also provide simplified proofs of many  results from~\cite{FY}.

A closed smooth manifold is  said
to be almost nonnegatively curved
if
it can Gromov--Hausdorff converge to a single point
under a lower curvature bound.
By rescaling, this definition is equivalent to the following one,
which we will employ thruout this article.

\begin{defn}
A closed smooth manifold $M$ is called almost nonnegatively curved
if it admits a sequence of Riemannian metrics $\{g_n\}_{n\in\mathbb{N}}$
whose sectional curvatures and diameters have the following bounds
$$\sec(M,g_n)\ge -1/n\quad\text{and}\quad\diam(M,g_n)\le 1/n.$$
\end{defn}

Almost nonnegatively curved manifolds
generalize almost flat as well as nonnegatively
curved manifolds.
The main source of examples
comes from a theorem of Fukaya and Yamaguchi in \cite{FY}.
It states that
if $F\to E\to B$ is a fiber bundle
over an almost nonnegatively curved manifold $B$
whose fiber $F$ is compact
and admits a nonnegatively curved metric 
which is invariant under the
structure group, 
then the total space $E$ is almost nonnegatively curved.
Further examples are given by closed manifolds which admit
cohomogeneity one actions of
compact Lie groups 
(compare~\cite{ST}).

In this work we combine collapsing techniques with a non-smooth analogue of the gradient flow of concave functions on Alexandrov spaces.
This notion is based  on the construction  of gradient curves of $\lambda$-concave functions used in ~\cite{PP} and bears many similarities to the Sharafutdinov retraction~\cite{Shar}.
Gradient flows on general metric spaces were considered in~\cite{grad-flow-book}.
The gradient flow on Alexandrov spaces plays a key role in the proofs of two of the three main results in this paper, and we believe that it should also prove useful for dealing with other problems related to collapsing under a lower curvature bound.

\subsection{} Let us first briefly recall some  previously known results:

Let $M=M^m$ be an almost nonnegatively curved $m$-manifold.
\begin{enumerate}[$\diamond$]
\item Gromov proved in ~\cite{G1}
that the minimal number of generators of the fundamental
group $\pi_1M$ of $M$ can be estimated by a constant $C_1(m)$
depending only on $m$,
and in ~\cite{G2} that the sum of Betti numbers of $M$ with
respect to any field of coefficients
does not exceed some uniform constant $C_2=C_2(m)$.
\item In \cite{Yam},
Yamaguchi showed that, up to a finite cover, the manifold $M$
fibers over a $b_1(M;\mathbb{R})$-dimensionsal torus
and that $M$ is diffeomorphic to a torus if $b_1(M;\mathbb{R})\z=m$. (Here $b_1(M;\mathbb{R})$ denotes the first Betti number with real coefficients.)
\item In \cite{FY},
Fukaya and Yamaguchi proved that 
\begin{enumerate}[$\circ$]
\item $\pi_1M$ is almost nilpotent;
that is, $\pi_1M$ contains a nilpotent subgroup of finite index
\item $\pi_1M$ is $C_3(m)$-solvable;
that is, contains a solvable subgroup of index at most $C_3(m)$
\end{enumerate}
\item If a closed manifold has  negative Yamabe constant,
then it cannot volume collapse with scalar curvature bounded from below
(see \cite{Sch,Lb}).  In particular, no such manifold
can be almost nonnegatively curved.
\item The $\hat{A}$-genus of a closed spin manifold $X$
of almost nonnegative Ricci curvature satisfies the inequality
$\hat{A}(X)\le 2^{\dim X/2}$ (\cite{G5}, \cite{Ga}).
\end{enumerate}

\subsection{}
Our first result concerns the relation between curvature bounds and the actions of the fundamental group on the higher homotopy groups.

Recall that an action by automorphisms of a group $G$ on an abelian group $V$
is called nilpotent if $V$
admits a finite sequence of $G$-invariant subgroups
 $$V=V_0\supset V_1\supset\ldots\supset V_k=0$$
such that the induced action of $G$ on $V_i/V_{i+1}$ is trivial for any $i$.
A connected CW-complex $X$ is called \emph{nilpotent}
if $\pi_1X$ is a nilpotent group that operates nilpotently
on $\pi_kX$ for every $k\ge 2$.

Nilpotent spaces play an important role in topology
since they enjoy some of the best
homotopy-theoretic properties of simply connected spaces,
like a Whitehead theorem or reasonable Postnikov towers.
Furthermore, unlike the category of simply connected spaces,
the category of nilpotent ones is closed under
many constructions such as the based loop space functor
or the formation of function spaces,
and group-theoretic functors, like localization and completion,
have topological extensions in this category.

\begin{main}[Nilpotency Theorem]\label{intro:main1}
Let $M$ be a closed almost nonnegatively curved manifold.
Then a finite cover of $M$ is a nilpotent space.
\end{main}

It would be interesting to know
whether  the order of this covering
can be estimated solely in terms of the dimension of $M$.

\begin{ex}
Let  $h\colon \mathbb{S}^3 \times \mathbb{S}^3\to \mathbb{S}^3 \times \mathbb{S}^3$ be defined by
$$h\colon (x,y)\mapsto (x{\cdot} y,y{\cdot} x{\cdot} y).$$
This map  is a diffeomorphism with inverse given by
$$h^{-1}\colon (u,v)\mapsto (u^2{\cdot} v^{-1},v{\cdot} u^{-1}).$$
The induced map $h_*$ on $\pi_3(\mathbb{S}^3\times\mathbb{S}^3)$ is given by the matrix
$A_h=\left(\begin{smallmatrix} 1&1\\
1&2
\end{smallmatrix}\right)$.
Notice that the eigenvalues of $A_h$
are different from $1$ in absolute value.
Let $M$ be the mapping cylinder of $h$.
Clearly, $M$ has the structure of a fiber bundle 
\[\mathbb{S}^3 \times \mathbb{S}^3 \to M \to \mathbb{S}^{1},\] 
and the
action of $\pi_1M\cong\mathbb{Z}$ on $\pi_3M\cong \mathbb{Z}^2$
is generated by  $A_h$.
In particular, $M$ is not a nilpotent space and hence, by Theorem~\ref{intro:main1},
it does not admit almost nonnegative curvature.
This fact doesn't follow from any previously known results.
\end{ex}

\subsection{}
Our next main result provides an affirmative answer to a conjecture of Fukaya and Yamaguchi~\cite[Conjecture 0.15]{FY}.

\begin{main}[$C$-Nilpotency Theorem for $\pi_1$]\label{intro:main2}

Let $M$ be an almost nonnegatively curved $m$-manifold.
Then $\pi_1M$ is $C(m)$-nilpotent;
that is,
$\pi_1M$ contains a nilpotent subgroup of index at most $C(m)$.

\end{main}
Notice  that Theorem ~\ref{intro:main2} is new even
for manifolds of nonnegative curvature.

\begin{ex} For any $C>0$ there exist
prime numbers $p>q>C$ and a finite group $G_{pq}$
of order $p q$
which is solvable but not nilpotent.
In particular, $G_{pq}$ does not contain any nilpotent subgroup
of index less than or equal to $C$.

Whereas none of the results mentioned so far
excludes 
$G_{pq}$
from being the fundamental group of some
almost nonnegatively curved $m$-manifold,
Theorem~\ref{intro:main2} shows that for $C>C(m)$ none of the groups $G_{pq}$
can be realized as the fundamental group of
such a manifold.
\end{ex}

\subsection{}
In  ~\cite{FY} Fukaya and Yamaguchi also conjectured
that a finite cover
of an almost nonnegatively Ricci curved manifold $M$
fibers over a nilmanifold with a fiber
which has nonnegative Ricci curvature and whose fundamental group is finite.
This conjecture
was later  refuted  by Anderson ~\cite{An}.

It is, on the other hand,  very natural to consider this conjecture
in the context of almost nonnegative {\sl sectional} curvature.
In fact, here 
Yamaguchi's fibration theorem (\cite{Yam})
and the results of \cite{FY} easily imply 
that a finite cover of an almost nonnegatively curved manifold
admits a map onto
 a nilmanifold whose homotopy fiber is a  simply connected closed manifold.

From mere topology, it is, however,
not clear whether this homotopy fibration
can actually  always be made into a genuine fiber bundle.
Our next result shows that
this is indeed true,
and that
for manifolds of almost nonnegative {sectional}  curvature
Fukaya's and Yamaguchi's original conjecture essentially does hold.

\begin{main}[Fibration Theorem]\label{intro:main3}

Let $M$ be an almost nonnegatively curved manifold.
Then a finite cover $\tilde{M}$ of $M$
is the total space of a fiber bundle
 $$
 F\to \tilde{M}\to N
 $$
over a nilmanifold $N$ with a simply connected fiber $F$.
Moreover, the fiber $F$
is  almost nonnegatively curved
in the generalized sense as defined below.

 \end{main}

\smallskip

 \begin{defn} \label{defn:ann}A closed smooth manifold $M$
is called almost nonnegatively curved in the generalized sense
if for some nonnegative integer $k$
there exists  a sequence of complete Riemannian metrics $g_n$ on
$M\times \mathbb{R}^k$
and points $p_n\in M\times \mathbb{R}^k$
such that
\begin{enumerate}
\item the sectional curvatures of the metric balls of radius $n$
around $p_n$ satisfy
$$\sec(B_n(p_n))\ge -1/n;$$
\item for $n\to\infty$
the pointed Riemannian manifolds $((M\times \mathbb{R}^k, g_n), p_n)$
converge in the pointed Gromov--Hausdorff distance to $(\mathbb{R}^k,0)$;
\item
the regular fibers over $0$
are diffeomorphic to $M$ for all large $n$.
\end{enumerate}
\end{defn}

Due to Yamaguchi's fibration theorem \cite{Yam},
manifolds which are almost nonnegatively curved in the generalized sense
play the same central role in collapsing under a lower curvature bound
as almost flat manifolds do
in the Cheeger--Fukaya--Gromov theory
of collapsing with bounded curvature  (see \cite{CFG}).

It is not known whether
all manifolds which are
almost nonnegatively curved 
in the generalized sense
are almost nonnegatively curved.
Clearly, if $k=0$, this definition reduces to the standard one.
Moreover, it is easy to see that all results of the present article,
as well as all results about
almost nonnegatively curved manifolds mentioned earlier
(except possibly for the ones concerning 
the $\hat A$-genus and Yamabe constant),
hold for manifolds which are almost nonnegatively curved
in the sense of
Definition~\ref{defn:ann}.

\smallskip

\subsection{}
Let us describe the structure of the remaining sections of this article.

In section~\ref{s:gradpush}, after providing some necessary background
from Alexandrov geometry,
we introduce the gradient flow of the square of a distance function.
It serves as one of
the main technical tools in the proofs of  theorem~\ref{intro:main1}
and  theorem~\ref{intro:main2}.

In section \ref{sec:nilp} we prove Theorem~\ref{intro:main1} by a direct application of the gradient flow technique.

In section \ref{sec:c-nilp} we  prove Theorem~\ref{intro:main2}.
The proof is also based on the gradient flow, but
is more involved and employs further technical tools such as
``limit fundamental groups'' of Alexandrov spaces.

In section \ref{sec:fib} we prove Theorem~\ref{intro:main3}.
This section is completely independent from the rest of the article.

In section 6 we discuss some further open questions related to our results.

\smallskip

{\it Acknowledgements.}
We would like to express our thanks to the following people
for helpful conversations during the preparation of this work:
I.~Belegradek,
V.~Gorbunov, E.~Formanek, I.~Kapovitch, A.~Lytchak, R.~Matveyev, D.~Robinson,
D.~Sullivan,  B.~Wilking, and Yu.~Zarkhin.

We are also especially grateful to Yunhui Wu for noticing a mistake
in the  proof of Theorem ~\ref{thm:cnilp} in the original version of
this article.

\section{ Alexandrov geometry and the gradient flow}\label{s:gradpush}

This section
provides necessary background in Alexandrov geometry.
The results of sections 2.1--2.3 are mostly  repeated from  \cite{PP}, \cite{Ptr2} and \cite{Ptr2007}.
The reader may consult \cite{BGP} for a general reference on Alexandrov spaces.

\subsection{$\lambda$-concave functions}

\begin{defn}(for a space without boundary) Let $A$ be an Alexandrov space without boundary.
A Lipschitz function $f\colon A \rightarrow \mathbb{R}$ is called
$\lambda$-\emph{concave} if for any unit speed minimizing geodesic $\gamma$
in $A$, the function
$$f\circ\gamma(t)-\lambda{\cdot}  t^2/2$$
is concave.
\end{defn}

If $A$ is an Alexandrov space with boundary, then its double $\tilde A$ is also an Alexandrov space (see \cite[5.2]{Per}).
Let $\text{p}\colon \tilde A\to A$ be the canonical map.
Given a function $f$ on $A$, set $\tilde f=f\circ \text{p}$.

\

\begin{defn}{\it (for a space with boundary)}  Let $A$ be an Alexandrov space with boundary. A Lipschitz function $f\colon A \rightarrow \mathbb{R}$ is called
$\lambda$-\emph{concave} if for any unit speed
minimizing geodesic $\gamma$ in $\tilde A$, the function
$$\tilde f\circ\gamma(t)-\lambda{\cdot}  t^2/2 $$
is concave.
\end{defn}

\begin{rmk}Notice that the restriction of a linear function on $\mathbb{R}^n$ to a ball is not $0$-concave in this sense.
\end{rmk}

\begin{rmk}
In the above definitions, the Lipschitz condition is only technical.
With some extra work,
all results of this section can be extended to continuous functions.
\end{rmk}

\subsection{Tangent cone and differential} Given a point $p$ in an Alexandrov space  $A$, we denote by $T_p=T_p(A)$ the tangent cone at $p$.

If $d$ denotes the metric of an Alexandrov space $A$,
let us denote by $\lambda{\cdot}  A$ the space $(A,\lambda{\cdot}  d)$.
Let $i_\lambda\colon  \lambda{\cdot}  A\to A$ be the canonical map.
The limit of $(\lambda{\cdot}  A,p)$ for $\lambda\to\infty$ is  the tangent cone $T_p$ at $p$ (see \cite[ 7.8.1]{BGP}).

\begin{defn} For any function $f\colon A \rightarrow \mathbb{R}$ the function $d_pf\colon T_p \rightarrow \mathbb{R}$ such that
$$d_pf=\lim_{\lambda\to\infty} \lambda{\cdot} (f\circ i_\lambda-f(p))$$
is called the differential of $f$ at $p$.
\end{defn}

It is easy to see that for a $\lambda$-concave function $f$
the differential $d_pf$ is defined everywhere, and that $d_pf$
is a $0$-concave function on the tangent cone $T_p$.

\begin{defn} Given a $\lambda$-concave function $f\colon A \rightarrow \mathbb{R}$, a point $p\in A$
is called critical point of $f$ if $d_pf\le 0$.
\end{defn}

\subsection{Gradient curves}
With a slight abuse of notation we will call elements of the tangent cone $T_p$ the ``tangent vectors'' at $p$.
The origin of $T_p$ plays the role of the zero
vector and is denoted by $o=o_p$.
For a tangent vector $v$ at $p$ we define its absolute value $|v|$
as the distance $|ov|$ in $T_p$.
For two tangent vectors $u$ and $v$ at $p$ we can define
their ``scalar product''
$$\langle u, v\rangle=(|u|^2+|v|^2-|uv|^2)/2=|u|{\cdot} |v|{\cdot} \cos\alpha,$$ where $\alpha=\angle uov$ in $T_p$.

For two points $p,q\in A$ we define $\log_pq$ to be a tangent vector $v$
at $p$
such that $|v|=|pq|$ and such that the
direction of $v$ coincides with a direction from $p$ to $q$
(if such a direction is not unique, we choose any one of them).
Given a curve $\gamma(t)$ in $A$, we denote by $\gamma^+(t)$
the right and by
$\gamma^-(t)$ the left tangent vectors to  $\gamma(t)$, where, respectively,
$$\gamma^\pm(t)\in T_{\gamma(t)},\quad\gamma^\pm(t)=\lim_{\eps\to +0}\frac{\log_{\gamma(t)}\gamma(t\pm\eps)}{\eps}.$$
For a real function $f(t)$, $t\in \mathbb{R}$, we denote by  $f^+(t)$ its right derivative and by $-f^-(t)$ its left derivative. Note that our sign convention (which is chosen to agree with the notion of right and left derivatives of curves) is not quite standard. For example,
$$\text{if}\quad f(t)=t\quad\text{then}\quad f^+(t)\equiv 1\quad\text{and}\quad f^-(t)\equiv -1.$$

\begin{defn}Given a $\lambda$-concave function $f$ on $A$,
a vector $g\in T_p(A)$ is called a gradient of $f$ at $p\in A$
(in short:  $g=\nabla_p f$) if

(i) $d_pf(x)\le \langle g , x\rangle$ for any $x\in T_p$, and

(ii) $d_pf(g) = \langle g,g \rangle .$
\end{defn}

It is easy to see that any $\lambda$-concave function
has a uniquely defined gradient vector field.
Moreover, if $d_pf(x) \le 0$ for all $x\in T_p$,
then $\nabla_p f=o$
(here $o$ denotes the origin of the tangent cone $T_p$); otherwise,
 $$\nabla_p f=d_pf(\xi){\cdot}  \xi$$
where $\xi$ is the (necessarily unique)
unit vector for which the function $d_pf$ attains its maximum.

Moreover, for any minimizing geodesic $\gamma\colon  [a,b]\to U$
parameterized by arclength, the following inequality holds:
\begin{equation}\label{e:star}
\langle \gamma^+(a),\nabla_{\gamma(a)} f\rangle +\langle\gamma^-(b),\nabla_{\gamma(b)} f\rangle
\ge -\lambda{\cdot} (b-a) .
\end{equation}

Indeed,
\begin{align*}
\langle   \gamma^+(a),  \nabla_{\gamma(a)}  f   \rangle +
\langle  \gamma^-(b),   \nabla_{\gamma(b)} f  \rangle
&\ge
d_{\gamma(a)}f(\gamma^+(a))+d_{\gamma(b)}f(\gamma^-(b))
=
\\
&=
(f\circ\gamma)^+|_a+(f\circ\gamma)^-|_b\ge
\\
&\ge -\lambda{\cdot} (b-a).
\end{align*}

\begin{defn}A curve $\alpha\colon [a,b] \rightarrow A$ is called an $f$-\emph{gradient curve} if for any $t\in [a,b]$
$$\alpha^+(t)=\nabla_{\alpha(t)}f.$$
\end{defn}

\begin{prop}\label{prop:gradlip} Given a $\lambda$-concave function $f\colon A \rightarrow \mathbb{R}$ and a point $p\in A$ there is a unique gradient curve $\alpha\colon [0,\infty) \rightarrow A$ such that $\alpha(0)=p$.

Moreover, if $\alpha$ and $\beta$ are two $f$-gradient curves, then
$$|\alpha(t_1)\beta(t_1)|\le |\alpha(t_0)\beta(t_0)|{\cdot} \exp(\lambda{\cdot}  (t_1-t_0))
\quad\text{for all}\quad
t_1\ge t_0.$$

\end{prop}

The gradient curve  can be constructed as a limit of broken geodesics,
made up of short segments with directions close to the gradient.
The convergence, uniqueness, as well as the last inequality
in Proposition~\ref{prop:gradlip}
follow from
inequality (\ref{e:star}) above, while Corollary~\ref{cor:gradlim}
below guarantees that the limit is indeed
a gradient curve, having a unique right tangent vector at each point.

\begin{lem}\label{lem:gradcon}
 Let $A_n \GHto A$ be a sequence of
Alexandrov spaces with curvature $\ge k$
which Gromov--Hausdorff converges to an Alexandrov space $A$.

Let $f_n\to f$, where $f_n\colon A_n\to \mathbb{R}$ is a sequence of $\lambda$-concave
functions converging to $f\colon A\to \mathbb{R}$.

Let $p_n\to p$,  where $p_n\in A_n$ and $p\in A$.

Then
$$|\nabla_p f|\le \liminf_{n\to \infty} |\nabla_{p_n} f_n|.$$
\end{lem}

\begin{cor}\label{cor:gradlim} 
For any $\lambda$-concave function $f$ on $A$
the function 
$$p\mapsto|\nabla_pf|$$ 
is lower semicontinuous;
that is, for any sequence of points $p_n\in A$, $p_n\to p$, we have
$$|\nabla_p f|\le \liminf_{n\to \infty} |\nabla_{p_n} f|.$$

\end{cor}

\begin{proof}[Proof of Lemma~\ref{lem:gradcon}] Fix an $\eps>0$ and choose $q$ near $p$ such
that
$$\frac{f(q)-f(p)}{|pq|}> |\nabla_p f|-\eps.$$
Now choose $q_n\in A_n$ such that $q_n\to q$.
If $|pq|$ is sufficiently small and $n$ is
sufficiently large,  the $\lambda$-concavity of $f_n$ then implies that
$$\liminf_{n\to \infty}\frac{d_{p_n}f_n(v_n)}{|v_n|}\ge |\nabla_p f|-2{\cdot} \eps
\quad\text{for}\quad 
v_n=\log_{p_n}(q_n)\in T_{p_n}(A_n).$$
Therefore,
$$\liminf_{n\to \infty} |\nabla_{p_n} f_n|\ge
|\nabla_p f|-2{\cdot} \eps
\quad\text{for any}\quad
\eps>0;$$
that is,
$$\liminf_{n\to \infty} |\nabla_{p_n} f_n|\ge
|\nabla_p f|. $$
\end{proof}

\begin{lem}\label{lem:concave}Let $f$ be a $\lambda$-concave function,  $\lambda\ge 0$ and $\alpha(t)$ be an $f$-gradient curve, and let $\bar\alpha(s)$ be its reparameterization by arclength. Then $f\circ\bar\alpha$ is  $\lambda$-concave.
\end{lem}

\begin{proof}

$$(f\circ\bar\alpha)^+(s_0)=|\nabla_{\bar\alpha(s_0)}f|
\ge
 \frac{d_{\bar\alpha(s_0)}f(\log_{\bar\alpha(s_0)}(\bar\alpha(s_1))}
 { |\bar\alpha(s_1)\,\bar\alpha(s_0)|}
 \ge$$
$$
\ge\frac{f(\bar\alpha(s_1))-f(\bar\alpha(s_0))
-
\lambda{\cdot}  |\bar\alpha(s_1)\,\bar\alpha(s_0)|^2/2 }{|\bar\alpha(s_1)\,\bar\alpha(s_0)|}
\ge $$
$$\ge
\frac{f(\bar\alpha(s_1))-f(\bar\alpha(s_0))}{ s_1-s_0}
-
\lambda{\cdot} |\bar\alpha(s_1)\,\bar\alpha(s_0)|/2.$$
Since $\frac{|\bar\alpha(s_1)\,\bar\alpha(s_0)|}{(s_1-s_0)}\to 1$
 as $s_1\to s_0+$, it follows that $f\circ\bar\alpha$ is  $\lambda$-concave.

\end{proof}

\subsection{ The Gradient Flow on Alexandrov Spaces}
Let $f$ be a $\lambda$-concave function on an Alexandrov space $A$.
Consider the map $\Phi_f^T\colon A\to A$ defined as follows:
$\Phi_f^T(x)=\alpha_x(T)$, where $\alpha_x\colon [0,\infty)\to A$ is the $f$-gradient curve with $\alpha_x(0)=x$.
The map $\Phi_f^T$ is called
\emph{$f$-gradient flow at time $T$}. From Proposition~\ref{prop:gradlip}  it is clear that  $\Phi_f^T$ is an $\exp(\lambda{\cdot}  T)$-Lipschitz map. Next we want to prove that this map behaves nicely under Gromov--Hausdorff-convergence.

\begin{thm}\label{thm:grconv}
Let $A_n \GHto A$
be a sequence of Alexandrov spaces with curvature $\ge k$
which converges to an Alexandrov space $A$.

Let $f_n\to f$, where $f_n\colon A_n\to \mathbb{R}$ is a sequence of $\lambda$-concave functions and $f\colon A\to \mathbb{R}$.

Then  $\Phi_{f_n}^T\to\Phi_f^T$.

\end{thm}

Theorem~\ref{thm:grconv}  immediately follows  from the following Lemma:

\begin{lem}Let $A_n \GHto A$
be a sequence of Alexandrov spaces with curvature $\ge k$
which converges to an Alexandrov space $A$.

Let $f_n\to f$, where $f_n\colon A_n\to \mathbb{R}$ is a sequence of $\lambda$-concave functions and $f\colon A\to \mathbb{R}$.

Let $\alpha_n\colon  [0,\infty) \to A_n$ be the sequence of
$f_n$-gradient curves with $\alpha_n(0)=p_n$ and
let $\alpha\colon  [0,\infty) \to A$ be the $f$-gradient curve with $\alpha(0)=p$.

Then $\alpha_n\to\alpha$.
\end{lem}

\begin{proof}
We may assume without loss of generality  that $f$ has no critical points.
(Otherwise consider instead the  sequence $A'_n=A_n\times\mathbb{R}$ with $f'_n(a\times x)=f_n(a)+x$.)

Let $\bar\alpha_n(s)$ denote the reparameterization of $\alpha_n(t)$
by arc length.
Since all $\bar\alpha_n$ are $1$-Lipschitz,
 we can choose a converging subsequence from
any subsequence of $\bar\alpha_n$.
Let $\bar\beta\colon [0,\infty)\to A$ be its limit.

Clearly, $\bar\beta$ is also 1-Lipschitz and hence $|\bar\beta^+|\le 1$.
Therefore, by Lemma \ref{lem:gradcon},
\begin{align*}
\lim_{n\to\infty}f_n\circ\bar\alpha_n|_a^b
&=
\lim_{n\to\infty}\int_a^b|\nabla_{\bar\alpha_n(s)} f_n|{\cdot}ds\ge
\\
&\ge\int_a^b|\nabla_{\bar\beta(s)} f|{\cdot}ds\ge
\\
&\ge
\int_a^b d_{\beta(t)} f(\beta^+(t)){\cdot} dt=f\circ\beta|_a^b.
\end{align*}

On the other hand, since $\bar\alpha_n\to\bar\beta$ and $f_n\to f$ we have
$f_n\circ\bar\alpha_n|_a^b \to f\circ\bar\beta|_a^b$.
Therefore, in both of these inequalities in fact equality holds.

Hence, $|\nabla_{\bar\beta(s)} f|= \lim_{n\to\infty} |\nabla_{\bar\alpha_n(s)} f_n|$,
$|\bar\beta^+(s)|= 1$
and the directions of $\bar\beta^+(s)$ and  $\nabla_{\bar\beta(s)} f$
coincide almost everywhere.
This implies that $\bar\beta(s)$ is a gradient curve reparameterized by
arc length.
In other words, if $\bar\alpha(s)$ denotes the
reparameterization of $\alpha(t)$
by arc length,
then $\bar\beta(s)=\bar\alpha(s)$ for all $s$.
It only remains to show that the original
parameter $t_n(s)$ of $\alpha_n$ converges to the original
parameter $t(s)$ of $\alpha$.

Notice that $|\nabla_{\bar\alpha_n(s)} f_n|{\cdot} dt_n=ds$ or
$dt_n/ds=ds/d(f_n\circ\bar\alpha_n)$.
Likewise, $dt/ds\z=ds/d(f\circ\bar\alpha)$.
Then the convergence $t_n\to t$ follows from the $\lambda$-concavity of
$f_n\circ\bar\alpha_n$ (see Lemma~\ref{lem:concave})
and the convergence $f_n\circ\bar\alpha_n\to f\circ\bar\alpha.$
\end{proof}

\subsection{Gradient balls}

\

Let $A$ be an Alexandrov space and let $S\subset A$ be a subset of $A$. A function $f\colon A\to \mathbb{R}$  which can be represented as
$$f=\sum_i \theta_i{\cdot} \frac{\dist^2_{a_i}}{2}
\quad\text{with }\quad
\theta_i\ge 0,\quad
\sum_i \theta_i=1
\quad\text{and}\quad
a_i\in S$$
will be called {\sl cocos-function with respect to $S$}
(where ``cocos'' stands for
{\bf co}nvex {\bf co}mbination of {\bf s}quares of distance functions).
A broken gradient curve for a collection of such functions
will be called cocos-curve with respect to $S$.

For $p\in A$ and $T,r\in \mathbb{R}_+$,
let us define  ``the gradient ball with center $p$
and radius $T$ with respect to $B_r(p)$'', $\beta^r_T(p)$,
as the set of all end points of cocos-curves with respect to $B_r(p)$
that start at $p$
with total time $\le T$.

\begin{lem}\label{lem:gradb}\

\begin{enumerate}[(I)]
\item There exists $T=T(m)\in \mathbb{R}_+$
such that for any $m$-dimensional Alexandrov space $A$ with curvature $\ge -1$ and any $q\in A$ there is a point $p\in A$
such that
\begin{enumerate}[(i)]
\item $|pq|\le 1$, and
\item $B_1(p)\subset \beta^{1}_T(p)$.
\end{enumerate}
\item There exists $T'=T'(m)\in \mathbb{R}_+$
such that the following holds.
Let $A$ be an Alexandrov space
which is a quotient $A=\tilde A/\Gamma$
of an  $m$-dimensional Alexandrov space $\tilde A$ with curvature $\ge -1$
by a discrete action of a group of isometries $\Gamma$.
Let $q\in A$ and $p=p(q)\in A$
be as in part I above.

Then
for any lift $\tilde p \in \tilde A$ of $p$ one has that
$B_1(\tilde p)\subset \beta^{1}_{T'}(\tilde p)$.
\end{enumerate}
\end{lem}

\begin{proof}
The proof is similar to the construction of a strained point in an  Alexandrov space (see~\cite{BGP}).

Set $\delta=10^{-m}$.
Take $a_1=q$ and take $b_1$ to be a farthest point from $a_1$ in the closed ball $\bar B_1(a_1)$.
Take $a_2$ to be a midpoint of $a_1b_1$ and let  $b_2$ be
a farthest point from $a_2$
such that $|a_1b_2|=|a_1a_2|$ and $|a_2b_2|\le\delta{\cdot} |a_1b_1|$, etc.
On the $k$-th step we have to take $a_{k}$ to be a midpoint of
$a_{k-1}b_{k-1}$
and $b_{k}$
to be a farthest point from $a_{k}$
such that $|a_ib_{k}|=|a_ia_{k}|$ for all $i< k$  and
$|a_{k}b_{k}|\le\delta{\cdot}  |a_{k-1}b_{k-1}|$.

After $m$ steps, take $p$ to be a midpoint of $a_mb_m$.
We only have to check that we can find $T=T(m)$
such that $\beta^{1}_T(p)\supset B_1(p)$.

Let $t_i$ be the minimal time such that $B_{|a_ib_i|/\delta^m}(p)\subset \beta^{1}_{t_i}(p)$.
Then one can take $T=t_1$.
Therefore it is enough to give estimates for $t_m$ and $t_{k-1}/t_k$ only in terms of $\delta$ and $m$.
Looking at the ends of broken gradient curves
starting at $p$
for the functions $\dist^2_{p}/2$, $\dist^2_{a_i}/2$
and $ \dist^2_{b_i}/2$, we easily see that $t_n\le 1/\delta^m$.
Now, looking
at the ends of broken gradient curves
starting at $ B_{|a_{k-1}b_{k-1}|/\delta^m}(p)$ for
the functions $\dist^2_{p}/2$, $\dist^2_{a_i}/2$
and $ \dist^2_{b_i}/2$,
we have that $t_{k-1}/t_{k}\le 1/\delta^m$.
Therefore $t_1\le 1/\delta^{m^2}=10^{-m^3}$.
This finishes the proof of part (I).

For part (II),
notice that
\begin{enumerate}[a)]
\item for any $r,t>0$ we have
$\beta^r_t(p)\subset B_{re^t}(p);$
\item if $\beta^r_t(p)\supset B_{\rho}(p)$, then
$\beta^r_{t+\tau}(p)\supset B_{\rho e^\tau}(p);$
\item if $\rho=|px|$ and $x\in \beta^{r+\rho}_t(p)$,
then $\beta^{r}_\tau(x)\subset \beta^{r+\rho}_{t+\tau}(p)$.
\end{enumerate}
Take $\eps=e^{-T}/4$ and apply part (I) of the lemma to $\frac1\eps{\cdot} \tilde A$
to find a point $p'\in \tilde A$ such that $|\tilde pp'|\le \eps$
and $B_{\eps}(p')\subset \beta^{\eps}_{T}(p')\subset \tilde A$.
Then for some deck transformation $\gamma$
we have $\gamma(p')\in \beta^\eps_T(\tilde p)\subset B_{\eps{\cdot}  e^T}(\tilde p)$.
Therefore it holds that $\gamma(p')\in B_{1/2}(\tilde p)$.
Hence, taking
$$T'=2{\cdot} T+1/\eps=2{\cdot} T+4{\cdot} e^{T},$$
we obtain
$$\beta^{1}_{T'}(\tilde p)\supset \beta^\eps_{T+1/\eps}(\gamma(p'))
\supset B_1(\tilde p).$$
\end{proof}

\subsection{Short basis.}\label{short-basis} We will use the following construction due to Gromov.

Given an Alexandrov space $A$ with a marked point $p\in A$,  and  a group $\Gamma$ acting discretely on $A=(A,d)$ one can define a short basis of the action of $\Gamma$ at $p$ as follows:

For $\gamma\in \Gamma$ define the norm of $\gamma$ by the formula
$|\gamma|= d(p,\gamma(p))$.
Choose $\gamma_1\in \Gamma$
with  the minimal norm in $\Gamma$.
Next choose $\gamma_2$ to have
minimal norm in $\Gamma\backslash \langle\gamma_1\rangle$.
On the $n$-th step choose
$\gamma_n$ to have minimal norm in
$\Gamma\backslash \langle\gamma_1,\gamma_2,\dots,\gamma_{n-1}\rangle$.
The sequence $\{\gamma_1,\gamma_2,\dots\}$
is called a \emph{short basis} of $\Gamma$ at $p$.
In general, the number of elements of a short basis can be finite or infinite.
In the special case of the action of the fundamental group $\pi_1(A,p)$  on the universal cover of $A$
 one speaks of the short basis of  $\pi_1(A,p)$.

It is easy to see that for a short basis
$\{\gamma_1,\gamma_2,\dots\}$ of the fundamental group of an Alexandrov space $A$
the following is true:
\begin{enumerate}
\item If $A$ has diameter $d$ then $|\gamma_i|\le 2{\cdot} d$.
\item If $A$ is compact then $\{\gamma_i\}$ is finite.
\item For any $i>j$ we have $|\gamma_i|\le |\gamma_j^{-1}\gamma_i|$.
\end{enumerate}
The third property implies that if
$\tilde p \in \tilde A$
is in the preimage of $p$ in the universal cover $\tilde A$ of $A$
and $\tilde p_i=\gamma_i(\tilde p)$, then
$$|\tilde p_i\tilde p_j|\ge
\max\{|\tilde p\tilde p_i|,|\tilde p\tilde p_j|\}.$$

As was observed by Gromov, if $A$ is an Alexandrov space with curvature $\ge \kappa$ and diameter $\le d$, the last inequality implies that $\angle \tilde p_i \tilde p \tilde p_j>\delta=\delta(\kappa,d)>0$.
This yields an upper bound on the number of elements
of a short basis in terms of $\kappa,d$ and the dimension of $A$.

\section{Nilpotency of almost nonnegatively curved manifolds}\label{sec:nilp}

In this section we prove Theorem~\ref{intro:main1}.

\subsection{Preliminary lemmas}
Let $M$ be an almost nonnegatively curved manifold.
Let us denote by $M_n=(M,g_n)$, $n\in \mathbb{N}$,
a sequence of Riemannian metrics on $M$ such that $\sec(M_n)\ge -1/n$
and $\diam(M_n)\le 1/n$.
Let $\tilde M\to M$ be the universal covering  and 
$\tilde M_n\to M_n$ be
the universal Riemannian covering of $M_n$
(that is, $\tilde M_n$ is
$\tilde M$ equipped with the pullback of the Riemannian metric $g_n$).

\begin{klem}\label{l:key} Given $\eps>0$ and $r_2>r_1>0$,
let
$\tilde M_n\supset B_{r_2}(p_n)\supset B_{r_1}(p_n)$.
Then, for $n$ sufficiently large,
there is a  $(1+\eps)$-Lipschitz map
\[\Phi_n\colon B_{r_2}(p_n)\to B_{r_1}(p_n)\]
which is homotopic to the identity on $B_{r_2}(p_n)$.
\end{klem}

\begin{proof} Fix $R>\!\!>r_2$
(here $R>1000(1+1/\eps)r_2$ will suffice).
Notice that as $n\to\infty$, we have that $B_R(p_n)\GHto B_R\subset \mathbb{R}^q$.
Choose a finite $R/1000$-net $\{a_i\}$ of $\partial B_R\z\subset \mathbb{R}^q$.
Let $a_{i,n}\in M_n$ be sequences such that $a_{i,n} \to a_n$.
Consider the sequence of functions
$f_n\colon M_n\to \mathbb{R}$ with $f_n=\min_i\{\text{dist}^2_{a_{i,n}}\}$.

For large $n$, the functions $f_n$ are $2$-concave in $B_R(p_n)$,
so that, in particular, the gradient flows
$\Phi_{f_n}^T|_{B_{r_2}(p_n)}$ are $e^{2{\cdot} T}$-Lipschitz.
Moreover, if $\xi_x$ denotes
the starting vector of a unit speed shortest geodesic from $x$ to $p_n$,
then  for any $x\in B_{r_2}(p_n) \backslash B_{r_1}(p_n)$
we have $\langle\xi_x,\nabla f \rangle\ge R/2$.
Therefore, if $T=2{\cdot} r_2/R$, then $\Phi_{f_n}^T(B_{r_2}(p_n) )\subset B_{r_1}(p_n)$.
Thus $\Phi_n=\Phi_{f_n}^{2{\cdot} r_2/R}$ provides a
$4r_2/R$-Lipschitz map $B_{r_2}(p_n)\to B_{r_1}(p_n)$,
and it is $(1+\eps)$-Lipschitz if one chooses $R$ sufficiently large.
\end{proof}

For $\gamma\in \pi_1M$, set $|\gamma|_n=|p\,\gamma(p)|_{\tilde M_n}$, see~\ref{short-basis}.

\begin{cor}\label{cor:norm} Let $M$ be almost nonnegatively curved manifold.
Let $$h\colon \pi_1M\to Aut(H^*(\tilde M,\mathbb{Z})/tor$$
be the natural  action of $\pi_1M$ on $H^*(\tilde M,\mathbb{Z})$.
Then there is a sequence of norms $|\!|*|\!|_n$ on
$H^*(\tilde M,\mathbb{Z})/tor$
such that the following holds.
Given any $\eps>0$, there is $n\in \mathbb{Z}_+$
such that for any  $\gamma\in \pi_1M$
with $|\gamma|_n\le 2{\cdot} \diam(M_n)$ we have  $|\!|h(\gamma)|\!|_n\le 1+\eps$.

\end{cor}

\begin{proof} \cite[theorem 0.1]{FY} and
Yamaguchi's fibration theorem ~\cite{Yam} imply
that if $n$ is sufficiently large,
for any fixed $r\in \mathbb{R}_+$
we have that for any $p_n\in \tilde M_n$
the inclusion map $B_r(p_n)\to \tilde M_n$ is a homotopy equivalence.

Let $|\!|*|\!|_{n,r}$ denote the $L_\infty$-norm
on differential forms on $B_{r}(p_n)\subset \tilde M_n$.

Fix $r_2>r_1>0$.
If $\omega$ is a differential form on  $B_{r_1}(p_n)\subset M_n$
and $n$ is sufficiently large,  Lemma~\ref{l:key} implies that
$$|\!|\Phi^*_n(\omega)|\!|_{n,r_2}
\le
(1+\eps){\cdot}|\!|\omega |\!|_{n,r_1}
\quad\text{and}\quad2{\cdot} \diam(M_n)\le r_2-r_1.$$

If now $\omega$ is a form on
$B_{r_2}(p_n)\in \tilde M_n$ and
$\gamma\in \pi_1M$
such that
$$|\gamma|_n=|p_n\,\gamma(p_n)|\le 2{\cdot} \diam(M_n)\le  r_2-r_1,$$
then $B_{r_1}(p_n)\subset B_{r_2}(\gamma(p_n))\subset \tilde M_n$, whence
$$|\!|\Phi^*_n(\gamma^*(\omega))|\!|_{n,r_2}\le
(1+\eps){\cdot}|\!|\gamma^*(\omega)|\!|_{n,r_1}\le
(1+\eps){\cdot}|\!|\omega|\!|_{n,r_2}.$$

Thus, for the induced norms on the de Rham cohomology
of $\tilde M$ (and on its integral subspace $H^*(\tilde M,\mathbb{Z})/tor)$)  we have
$$|\!|[\gamma^*(\omega)]|\!|_{n,r_2}\le
(1+\eps){\cdot}|\!|[\omega]|\!|_{n,r_2}.$$

Therefore the sequence of norms $|\!|*|\!|_n=|\!|*|\!|_{n,r_2}$ satisfies the conditions of the Corollary.
\end{proof}

\begin{lem}\label{lem:eigen}
There exists a constant $N=N(n,k)\in \mathbb{Z}_+$ such that the following holds.
If $G$ is a subgroup of $\mathrm{GL}(n,\mathbb{Z})$ and $S$ is a set of generators of $G$
with $\#(S)\le k$ such that the eigenvalues of each element of
$S^N$ are all equal to $1$ in absolute value,
then the same is true for the eigenvalues of all
elements of $G$.
\end{lem}

\begin{proof}

Let $B$ be the set of all matrices in $\mathrm{GL}(n,\mathbb{Z})$
for which all of their eigenvalues are equal to 1 in absolute value.
Since the characteristic polynomials of such matrices
are uniformly bounded and have integer coefficients,
there are only finitely many of them.
Let $\bar{B}$ be the Zariski closure of $B$ in the set of all real
$n\times n$ matrices. By the above,
all elements of $\bar{B}$ satisfy that the absolute values of all
of their eigenvalues are equal to 1.

Consider now the space $V=\mathbb{R}^{k{\cdot}n^{2}}$ of
$k$-tuples of real $n\times n$ matrices.

Consider a  collection of matrices $(M_1,M_2,\dots,M_k)\in V$,
where $M_i\in \mathrm{GL}(n,\mathbb{R})$.
Let $F_k$ be a free group on $k$ generators,
generated by $S=\{\gamma_1,\gamma_2,\dots,\gamma_k\}$,
and let $h\colon F_k\to \mathrm{GL}(n,\mathbb{R})$ be the
homomorphism defined by $h(\gamma_i)=M_i$.
The property that for any $\gamma\in F_k$ $h(\gamma)$
be an element of $\bar{B} $
then describes an algebraic subset $A_\gamma\subset V$.

The  intersection $A=\cap_{\gamma\in F_k}A_\gamma$
is also  algebraic, and therefore there is a finite number $N=N(n,k)$
such that for $S^N\subset  F_k$,
$A=\cap_{\gamma\in S^N}A_\gamma$.
\end{proof}

\begin{lem}\label{lem:eigen1}
Let $\Gamma$  be a subgroup of $\mathrm{GL}(n,\mathbb{Z})$ such that
the eigenvalues of each element of $\Gamma$
are equal to $1$ in  absolute value.
Then $\Gamma$ contains a   subgroup  $\Gamma'$ of finite index
all of whose  elements have eigenvalues equal to $1$.
\end{lem}

\begin{proof}
Let $G$ denote the Zariski closure of $\Gamma$ in  $\mathrm{GL}(n,\mathbb{R})$ .
Then $G$, being an algebraic group,
is a Lie group with finitely many components.
Let $G_\circ$ be the identity  component of $G$.
By the same argument as in the proof of the previous lemma,
the set of all characteristic polynomials
of the elements of $G$ is finite.
Therefore the characteristic polynomial
of any element of $G_\circ$
is identically equal to $(x-1)^n$.

Therefore, the subgroup $\Gamma'=\Gamma\cap G_\circ$
satisfies all conditions of the Lemma.
\end{proof}

\begin{rmk}
As was pointed out to us by Yu.~Zarkhin, one can alternatively take
$\Gamma'$ to be the kernel of the composition of the homomorphisms
$\Gamma\to \mathrm{GL}(n,\mathbb{Z})\z\to \mathrm{GL}(n,\mathbb{Z}/3\mathbb{Z})$.
In this way one obtains a bound
$$[\Gamma:\Gamma']\le 3^{n^2}.$$
To see that $\Gamma'$ satisfies the conclusion of
Lemma~\ref{lem:eigen1}, one should notice that every element of
$\Gamma$ is a quasi-unipotent matrix since all its eigenvalues are
roots of unity.
The desired result then follows from the so-called
Minkowski Lemma.
Apply, for instance, \cite[Th. 7.2]{Zar} for $n=3,\
k=1$ (so $R(1,3)=1$), where we take $\mathcal O$ to  be the ring of
$n\times n$ integer matrices.
\end{rmk}

\subsection{Proof of Theorem~\ref{intro:main1}}
Let $M$ be an almost nonnegatively curved manifold.
Denote, as usual, by $M_n=(M,g_n)$, $n\in \mathbb{N}$,
a sequence of Riemannian metrics on $M$ such that $\sec(M_n)\ge -1/n$
and $\diam(M_n)\le 1/n$,
by $\tilde M$ the universal covering of $M$, and by
$\tilde M_n\to M_n$
the universal Riemannian covering of $M_n$.

After passing to a finite cover of $M$, by~\cite{FY}
we may assume that
 $\pi_1M$ is nilpotent.

Fix $p\in M$ and
let $\{\gamma_{i,n}\}$ be a short basis of $\pi_1(M_n,p)$ (see~\ref{short-basis}).
Then, if $n$ is sufficiently large,
the short basis $\{\gamma_{i,n}\}$
has at most $k=k(\dim M)$ elements and
its elements satisfy
$|\gamma_{i,n}|_n\le 2/n$ for every $i$.
Moreover,
Corollary~\ref{cor:norm}  implies
that given $\eps>0$, for all large $n$ and every $i$
we have $|\!|h(\gamma_{i,n})|\!|_n <
1+\eps$ and $|\!|h(\gamma^{-1}_{i,n})|\!|_n<1+\eps$.

Take $N=N(k,m)$ as in Lemma~\ref{lem:eigen}.
One can choose $\eps>0$
so small that
if $p$ is a polynomial with integer coefficients
for which all of its roots have  absolute
values lying between  $1/(1+\eps)^N$ and $(1+\eps)^N$,
then  all roots of $p$ have absolute values equal to $1$.
This follows from the fact
that the total number of integer polynomials
all of whose roots are contained in a fixed bounded region is finite.

Set $S_n:=\{\gamma_{i,n}\}$.
Then for any $\gamma\in S_n^N$ we have
$|\!|h(\gamma)|\!|_n<(1+\eps)^N$ and
$|\!|h(\gamma^{-1})|\!|_n<(1+\eps)^N$.
Therefore the absolute values of all eigenvalues
lie between $1/(1+\eps)^N$ and $(1+\eps)^N$.
Since the characteristic polynomial of $h(\gamma)$
has integer coefficients,
the absolute values of all the eigenvalues of $h(\gamma)$ are
in fact equal to $1$.

Apply now Lemma~\ref{lem:eigen}. It follows
that for any $\gamma\in \pi_1M$
the absolute values of all  eigenvalues of $h(\gamma)$ are equal to $1$.

Then Lemma~\ref{lem:eigen1}  implies that after passing to a finite
cover $M'$ of $M$,
for any $\gamma\in \pi_1M'$
all eigenvalues of $h(\gamma)$ are equal to $1$.
By Engel's theorem, one can choose an integral basis
of $H^*(\tilde M,\mathbb{R} )$
such that the action of $\pi_1M$
on $H^*(\tilde M, \mathbb{Z})/tor$ is given by upper triangular matrices.

Therefore, by passing to a finite cover $M''$ of $M'$,
we can assume that the action of $\pi_1M''$
on $H^*(\tilde M,\mathbb{Z} )$  (and  on $H_*(\tilde M,\mathbb{Z} )$)  is nilpotent.

Recall (see, e.g., \cite[2.19]{HMR})
that a connected CW complex with nilpotent fundamental group
is nilpotent if and only if the action of
its fundamental group on the homology of its universal cover is nilpotent.
Thus $M''$ is a nilpotent space, whence the proof of Theorem~$A$
is complete.\qed

\section{C-nilpotency of the fundamental group}\label{sec:c-nilp}

\subsection{}In this section we will prove Theorem~\ref{intro:main2}. It will follow
from the following somewhat stronger result.

\begin{thm}\label{thm:cnilp}
For any integer $m$ there exist constants $\eps(m)>0$  and $C(m)>0$ such that
the following holds.
If
$M^m$ is a closed smooth $m$-manifold which admits
a Riemannian metric $g$ with
$\sec(M^m,g)> -\eps(m)$ and $\diam(M^m,g)< 1$,
then the fundamental group of $M^m$ is $C(m)$-nilpotent; 
that is,
$\pi_1M^m$ contains a nilpotent subgroup of index $\le C(m)$.
\end{thm}
\begin{rmk} 
The proofs of Theorems~\ref{intro:main1} and ~\ref{intro:main3} show
that
corresponding versions of those results
also do hold when these theorems are reformulated
in a fashion similar to Theorem~\ref{thm:cnilp}.
\end{rmk}
By an argument by contradiction,
Theorem~\ref{thm:cnilp} follows
from the following  statement:

Given a sequence of Riemannian $m$-manifolds $(M_n,g_n)$
such that
\[\diam(M_n,g_n)\le 1/n
 \quad\text{and}\quad
 \sec(g_n)\ge -1/n,
\]
for each $n$, one can find $C\in \mathbb{R}$  such that
$\pi_1M_{n}$ is $C$-nilpotent for all sufficiently large $n$.

\subsection{Algebraic lemmas}

Recall that
the group of outer automorphisms $\Out(G)$ of a group $G$
is defined as the quotient of its automorphism group
$\Aut(G)$ by the subgroup of inner automorphisms $\Inn(G)$.

\begin{lem}[\it A characterization of $C$-nilpotent groups]\label{lem:cnilp}
Let
$$\{1\}=G_\ell\subseteq \ldots\subseteq G_1\subseteq G_0=G$$
  be a sequence of groups satisfying the following properties:

For any $i$ we have that
\begin{enumerate}[(i)]
\item $G_i \unlhd G$ is normal in $G$;
\item the image of the conjugation homomorphism
$h_i\colon  G\to \Out(G_i/G_{i+1})$ is finite of order at most $C_i$;
\item $G_i/G_{i+1}$ contains an abelian subgroup $E_i$ of index  $\le c_i$.
\end{enumerate}
Then $G$ contains a nilpotent subgroup $N$ of index at most
$$C=C(c_1,\ldots, c_\ell, C_1,\ldots, C_\ell)$$
where $N$ is of nilpotency class $\le \ell$.
\end{lem}
\begin{proof}
First of all, notice that
property (i) assures that the
objects described in parts (ii) and (iii) of the lemma are well-defined.

Set $\Gamma_i:=G_i/G_{i+1}$.

Let $H=\cap\ker{h_i}$.
Notice that $[G{:}H]\le \prod_iC_i$ and that the image of $H$ under the conjugation homomorphism $f_i\colon  G\to \Aut(\Gamma_i)$  lies in $\Inn(\Gamma_i)$;
that is,  $f_i|_H\colon H\z\to\Inn(\Gamma_i)$.

By passing to a subgroup we can assume that $E_i\unlhd \Gamma_i$ is normal of index $\le C(c_i)$ (we can take $C(c_i)=c_i!$).

By increasing $E_i$ if necessary we can assume that $E_i$ contains the center of $\Gamma_i$.

Let $Z_i$ be the image of $E_i$ under the projection map $\pi\colon \Gamma_i\to \Inn(\Gamma_i)$.
Clearly $[ \Inn(\Gamma_i): Z_i]\le c_i$ and $Z_i\unlhd  \Inn(\Gamma_i)$ is normal.

Let
$$N=H\cap\left(\bigcap_if_i^{-1}(Z_i)\right),\quad N_i=N\cap G_i.$$
Then
$$[G:N]\le C=C(c_1,\ldots, c_\ell, C_1,\ldots, C_\ell)$$
and $N$ satisfies:

For any $i$
\begin{enumerate}[(i$'$)]
\item$N_i \unlhd N$ is normal in $N$;

\item $N_i/N_{i+1}$ is in the center of $N/N_{i+1}$;
\end{enumerate}
 that is, $N$ is nilpotent of nilpotency length $\le \ell$.

 Condition (i$'$) is obvious so we only need to check (ii$'$).

 To see (ii$'$) observe that by construction the image of the conjugation action $N\to \Aut(\Gamma_i)$ lies in $\Inn(\Gamma_i)$. Moreover, in fact it lies in $\pi(A_i)$ and as such it acts trivially on $E_i$. Lastly observe that  $N_i/N_{i+1}\subset E_i$.

 Indeed, by construction, for any $g\in N_i/N_{i+1}\subset \Gamma_i$ there is $a\in E_i$ such that $\pi(g)=\pi(a)$. Therefore, $g=a{\cdot} z$ for some $z$ in the center of $\Gamma_i$. By our assumption on $E_i$ this means that $g\in E_i$.

 Thus $N$ acts trivially on $N_i/N_{i+1}$ which means that $N$ is nilpotent and $G$ is $C$-nilpotent.
 \end{proof}

\begin{tlem}[\it A characterization of finite actions]\label{lem:triv}
If $S$ is a finite set of generators of a group $G$ with $S^{-1}=S$,
and $h\colon G\to H$ is a homomorphism with $|h(S^n)|< n$ for some $n>0$,
then $h(S^n)=h(G)$ and, in particular, $|h(G)|< n$.
\end{tlem}

Let now $\Gamma$ be a group which acts discretely
by isometries on an Alexandrov space $A$ with curvature $\ge -1$.
Choose a marked point $p\in A$.
Assume that $\{\gamma_1,\gamma_2,\ldots,\gamma_n\}$
is a finite short basis of $\Gamma$ at $p$ (see~\ref{short-basis}),
and that $\theta\le|\gamma_i|\le 1$,
where $|\gamma|\overset{\rm def}{=}|p\gamma(p)|$.
Let $\#(R)$
denote the number of elements $\gamma\in\Gamma$ with $|\gamma|\le R$.
The  Bishop--Gromov inequality implies that
$$\#(R)\le v^m_{-1}(R)/v^m_{-1}(\theta),$$
where $m=\dim A$ and $v^m_{-1}(r)$
is the volume of the ball of radius $r$
in the $m$-dimensional simply connected space form of curvature $-1$.
Therefore,
if $\#\#(L)$ denotes the number of homomorphisms
$h\colon \Gamma\to \Gamma$ with norm $\le L$
(that is, the number of homomorphisms
for which it holds that for any $\gamma\in\Gamma$
one has that $|h(\gamma)|\le L|\gamma|$),
then
\begin{equation}\label{e:hombd}
\#\#(L)\le \#(L)^n\le \left[\frac{v^m_{-1}(L)}{v^m_{-1}(\theta)}\right]^n.
\end{equation}

\subsection{The blow-up construction}\label{sec:constr}

For $n\to\infty$, the manifolds $M_n$ clearly converge to a point $=:A_0$.

Set $M_{n,1}:=M_n$ and $\vartheta_{n,1}:=\diam M_{n,1}$.

Rescale now $M_{n,1}$ by $\tfrac{1}{\vartheta_{n,1}}$
so that $\diam(\tfrac{1}{\vartheta_{n,1}}\cdot M_{n,1})=1$.
Passing to a subsequence if necessary, one has
that the manifolds $\frac{1}{\vartheta_{n,1}}{\cdot} M_{n,1}$
converge to $A_1$,
where $A_1$
is a compact nonnegatively curved Alexandrov space with diameter $1$.

Now choose a regular point $\bar p_1\in A_1$,
and consider distance coordinates around $\bar p_1\in U_1\to \mathbb{R}^{k_1}$,
where $k_1$ is the dimension of $A_1$.
The distance functions can be lifted
to $U_{n,1}\subset \frac{1}{\vartheta_{n,1}}{\cdot} M_{n,1}$.

Let  $M_{n,2}$  be the level set of $U_{n,1}\to \mathbb{R}^{k_1}$
that corresponds to $\bar p_1$.
Clearly, $M_{n,2}$ is a compact submanifold of codimension $k_1$.
Set $\vartheta_{n,2}:=\diam M_{n,2}$.

Passing again to a subsequence if necessary, one has
that the sequence $\frac{1}{\vartheta_{n,2}}{\cdot} M_{n,2}$
converges to some  Alexandrov space $A_2$.
As before, $A_2$ is a compact nonnegatively curved Alexandrov space
with diameter 1. Set $k_2:=k_1+\dim A_2$.
If one now chooses a marked point in $M_{n,2}$,
then, as $n\to\infty$,
$M_n/\vartheta_{n,2}$ converges to $A_2\times\mathbb{R}^{k_1}$,
which is of some dimension $k_2>k_1$.

We repeat this procedure until, at some step, $k_\ell=m$.

\medskip

As a result one obtains
a sequence $\{A_i\}$ of compact nonnegatively curved Alexandrov spaces
with diameter $1$ that satisfies
$$\dim A_i=k_i-k_{i-1}, 
\quad\text{so that}\quad
\sum_{i=1}^\ell\dim A_i=m.$$
We also obtain
a sequence of rescaling factors $\vartheta_{n,i}=\diam M_{n,i}$,
and  a nested sequence of submanifolds
$$\{p_n\}=M_{n,\ell}\subset \cdots \subset M_{n,2}\subset M_{n,1}=M_n,$$
which in turn induces a sequence of homomorphisms
$$\{1\}=\pi_1M_{n,\ell}\buildrel {\imath}
\over \rightarrow \cdots\buildrel {\imath}
\over \rightarrow \pi_1M_{n,2}\buildrel {\imath}
\over \rightarrow \pi_1M_{n,1}=\pi_1M_n.$$

Let $G_i:=G_i(n):=\imath^i[\pi_1M_{n,i}]$.

For $n$ sufficiently large,
the subgroups $G_i(n)$ are those which are
generated by elements of  norm $\le 3{\cdot} \vartheta_{n,i}$.
Equivalently, if one takes a short basis $\{\gamma_i\}$ of $G(n)$,
then $G_i$ is the subgroup generated by all elements $\gamma_i$
such that $|\gamma_i|\le 3{\cdot} \vartheta_{n,i}$.

\subsection{ Limit fundamental groups of  Alexandrov spaces.}\label{fun-gr}

We will now define the ``limit'' or ``L-fundamental groups'' of
the Alexandrov spaces $A_i$ constructed above.
This notion is similar to the notion of  the fundamental group of an orbifold.
However, we note in advance
that the construction of the L-fundamental group does not only
depend on the spaces $A_i$,
but also on the chosen rescaled subsequence of $M_n$.
In fact, the following construction shows
that the limit fundamental group of $A_i$,
$\pi^L_1(A_i)$,
is isomorphic to $\pi_1(M_{n,i},M_{n,i+1})$
for all sufficiently large $n$.
But, unlike $\pi_1M_{n,i}$, the groups $\pi_1^LA_i$ will not depend on $n$.

{\it The limit fundamental groups of $A_i$.} 
Consider the converging sequence
\[(M_{n}/\vartheta_{n,i},p_{n})\GHto (A_i\times \mathbb{R}^{k_{i-1}},\bar p_i\times 0)\]
(here the interesting case is collapsing).
Recall that
$\bar p_i\in A_i$
is a regular point.
Fix $\eps>0$ such that
$\text{dist}_{\bar p_i}$ on $A_i$
does not have critical values in $(0,2{\cdot} \eps)$.
Take a sequence $R_n$ which converges very slowly to infinity
(here we will need $R_n\z\cdot \vartheta_{n,i}/\vartheta_{n,i-1}\to 0$ and $R_n\to \infty$).

Consider then
a sequence of Riemannian coverings
$$\Pi\colon (\tilde B_n,\tilde p_n)\to (B_{R_n}(p_n), p_n)$$
of
$B_{R_n}(p_n)\subset {M}_n/\vartheta_{n,i}$ with
$\pi_1(\tilde B_n,\tilde p_n)=\pi_1(B_\eps(p_{n}), p_n)$,
where $B_\eps(p_{n})\subset M_n/\vartheta_{n,i}$.

After passing to a subsequence if necessary,
the sequence $(\tilde B_n,\tilde p_n)$
converges to a nonnegatively curved Alexandrov space
$\tilde A_i\times \mathbb{R}^{k_{i-1}}$,
where the space $\tilde A_i$ has the same dimension as $A_i$.
Indeed, by construction it  contains an isometric copy of
$B_\eps(p_{n,i})$, and therefore
$$\dim \tilde A_i +k_{i-1}=
\dim \lim_{i\to\infty} B_\eps(p_{n,i})=\dim A_i+k_{i-1}.$$

Let us show that for all sufficiently large $n$,
$$\imath[\pi_1M_{n,i+1}]\trianglelefteq\pi_1M_{n,i}.$$
Assume that $\Pi(\tilde q_n)=\tilde p_n$
and that $\tilde q_n\to \bar q_n\in \tilde A_i$.
Connect $\bar p_n$ and $\bar q_n$ by a geodesic
which, by \cite{Ptr3}, only passes regular points.
Note that in a small neighborhood  of this geodesic in $M_n$
 we have two copies of $M_{n,i+1}$, near $\tilde p_n$ and $\tilde q_n$.
Therefore,
applying Yamaguchi's  Fibration Theorem
in a small neighborhood  of this geodesic,
we can construct a diffeomorphism from $M_{n,i+1}$ to itself.
This  implies that for any loop $\gamma$ which after lifting
connects $\tilde p_n \tilde q_n$,
we have 
$$\gamma^{-1}{\cdot} \imath[\pi_1M_{n,i+1}]{\cdot}
\gamma\subset \imath[\pi_1M_{n,i+1}];$$  
that is,
$$\imath[\pi_1M_{n,i+1}]\lhd\pi_1M_{n,i}$$
(for an alternative argument see also ~\cite{FY}).

This easily yields that $A_i=\tilde A_i/\Gamma_i$,
where $\Gamma_i$ is a group of isometries
which acts 
discretely on $\tilde A_i$.
The group $\Gamma_i$ is denoted  by $\pi_1^LA_i$
(the \emph{limit} or \emph{L-fundamental group} of $A_i$).
This group  is clearly isomorphic to
$$\pi_1(M_{n,i},M_{n,i+1})=\pi_1M_{n,i}/\imath[ \pi_1M_{n,i+1}]$$
for all sufficiently large $n$, and
the space $\tilde A_i$ will be called the \emph{universal covering} of $A_i$.

Since $\tilde A_i$ is nonnegatively curved
and $A_i=\tilde A_i/\pi_1^LA_i$ is compact,
by Toponogov's splitting theorem
$\tilde A_i$ isometrically splits
as $\tilde A_i=K_i\times \mathbb{R}^{s_i}$,
where $K_i$ is a compact Alexandrov space with curv$\ge 0$.
Since $\pi_1^LA_i$ is a group of isometries
that acts 
discretely on $\tilde A_i$,
it follows that $\pi_1^LA_i$ is a virtually abelian group.

\subsection{Final steps}

Consider now the corresponding series
$$\{1\}=G_\ell(n)\subset\ldots\subset G_1(n)\subset G_0(n)=\pi_1M_n.$$

The theorem then follows from the following

\

\begin{lem}\label{l:claim}
For all sufficiently large $n$, the series
$$\{1\}=G_\ell(n)\subset\ldots\subset G_1(n)\subset G_0(n)$$
constructed above
satisfies the assumptions of Lemma~\ref{lem:cnilp}
for numbers $C_i$ and $c_i$ which do not depend on $n$.
\end{lem}

We first prove the following
\begin{slem}\label{lem:normal}
 Each subgroup $G_i(n)$ is normal in $G(n)$.
\end{slem}

\begin{proof}
We will show by reverse induction on $k$
that $G_i(n)\unlhd G_{k}(n)$ for any $k\le i$.
Let us assume that we already know that $G_i(n)\unlhd G_{k+1}(n)$.
Since 
$$\imath[\pi_1M_{n,k+1}]\unlhd\pi_1M_{n,k},$$
we know that $G_{k+1}(n)\unlhd G_{k}(n)$.
Consider the covering
$$\Pi_{k+1}\colon (\tilde M_{n,k+1},\tilde p_{n,k+1})\to (M_n,p_n)$$
with covering group $\Gamma_{k+1}(n)$.

Clearly $(\tilde M_{n,k+1},\tilde p_{n,k+1})\GHto \mathbb{R}^{s_i}$
for some integer $s_i$.
Applying Lemma~\ref{lem:gradb}, it follows
that for any $a\in G$ with $|a|<1$ there is
a cocos-curve $\gamma$ in $\tilde M_{n,k+1}$
with total time $T$ connecting
$\tilde p_n$ and $a(\tilde p_n)$ in $\tilde M_{n,k+1}$.
Then clearly $\gamma\sim g{\cdot} a$ for some $g\in G_{k+1}(n)$.
Let us denote by $\Phi^T\colon  \tilde M_{n,i}\to \tilde M_{n,i}$ the gradient flow
corresponding  to $\gamma$.

Let $\gamma_j$ be a loop
from the short basis of $G_i(n)$.
As was mentioned in ~\ref{sec:constr},
if $n$ is large, then $\length\gamma_j\le 3{\cdot} \vartheta_{n,i}$.
Let us denote by  $\tilde\gamma_j$ a lift of $\gamma_j$ to $\tilde M_{n,i}$.
Let $\tilde p_{n,j}\in \tilde M_{n,i}$ be its starting point.
Since $[\gamma_j]\in G_i(n)$, we have that
$\tilde\gamma_j$ is a loop in $\tilde M_{n,i}$.
Consider then the loop
 $\gamma_j'=\Pi\circ\Phi^T\circ
\tilde\gamma_j$.
Clearly,
$$[\gamma_j]=a^{-1}{\cdot} g^{-1}{\cdot} [\gamma_j']{\cdot} g{\cdot} a ,
\quad\text{or}\quad[\gamma_j']=g{\cdot} a{\cdot} [\gamma_j]{\cdot} a^{-1}{\cdot} g^{-1}.$$

Now Proposition~\ref{prop:gradlip} implies that
$$\length(\gamma_j')\le \exp(2{\cdot} T){\cdot} \length(\gamma_j).$$
Thus, for sufficiently large $n$,
$$ g{\cdot} a{\cdot} [\gamma_j]{\cdot} a^{-1}{\cdot} g^{-1} \in G_i(n),$$

and since $g\in G_i(n)\unlhd  G_{k+1}(n)$ it follows that
$$ a{\cdot} [\gamma_j]{\cdot} a^{-1} \in G_i(n);$$
that is, $G_i(n)\unlhd  G_{k}(n)$.

\end{proof}

\begin{proof} [Proof of Lemma~\ref{l:claim}]

The group
$$\pi_1^LA_i=\pi_1(M_{n,i},M_{n,i+1})=\pi_1M_{n,i}/\imath[\pi_1M_{n,i}]$$
is virtually abelian.
Let $d_i$ be the minimal index of an abelian subgroup of $\pi_1^LA_i$.
The epimorphism
$\imath^i\colon \pi_1M_{n,i}\to G_i$
induces an  epimorphism $\pi_1^LA_i\z\to G_i(n)/G_{i+1}(n)$.
Therefore, $G_i(n)/G_{i+1}(n)$ is $d_i$-abelian for all large $n$.

Consider the covering
$\Pi_i\colon \tilde M_{n,i}\to M_n$ with covering group $G_i(n)$,
and let $\tilde p_{n,i}$ be a preimage of $p_n$.
Clearly $(\tilde M_{n,i},\tilde p_{n,i})\GHto \mathbb{R}^{s_i}$
for some integer $s_i$.
Applying Lemma~\ref{lem:gradb},
it follows  that for any $a\in G(n)$ with $|a|<1$ there is
a cocos-curve $\gamma$ in $\tilde M_{n,i}$ which connects $p$ and $a(p)$.
Then clearly $\gamma\sim g{\cdot} a$ for some $g\in G_i(n)$.
Let us denote by $\Phi^T\colon  \tilde M_{n,i}\to \tilde M_{n,i}$
the gradient flow corresponding  to $\gamma$.

Let $b\in G_i(n)$ and $\beta$ be a loop  representing $b$.
Let us denote by  $\tilde\beta$ a lift of $\beta$ to $\tilde M_{n,i}$.
Let $\tilde p_{n,i}\in \tilde M_{n,i}$ be its starting point.
Since $[\beta]\in G_i(n)$,
we have that $\tilde\beta$ is a loop in $\tilde M_{n,i}$.

Consider now the loop
$\beta'=\Pi\circ\Phi^T\circ\tilde\beta$. Clearly,
$$b=[\beta]=a^{-1}{\cdot} g^{-1}{\cdot} [\beta']{\cdot} g{\cdot} a , \quad\text{or}\quad[\beta']=g{\cdot} a{\cdot} b{\cdot} a^{-1}{\cdot} g^{-1}.$$
Proposition~\ref{prop:gradlip}  then implies that
\begin{equation*}\label{e:1}\length(\beta')\le \exp(2{\cdot} T){\cdot}\length(\beta)
\end{equation*}
Therefore, if $h_a\colon G_i(n)/G_{i+1}(n)\to G_i(n)/G_{i+1}(n)$
is induced by the conjugation $b\to a{\cdot} b{\cdot} a ^{-1}$,
then for any $a\in G(n)$ there is $g\in G_i(n)$ such that $|h_{g{\cdot} a}|\le \exp(2{\cdot} T)$.

Let now $\delta_i$ be the minimal norm of the elements of $\pi_1^LA_i$,
where $\pi_1^LA_i$ acts on $\tilde A_i$.
Then (\ref{e:hombd})  implies that the image of the
action of $G(n)$ by conjugation in $Out(G_i(n)/G_{i+1}(n))$ is $C_i$-finite,
where  $C_i$ depends only on $c_i$, $T$,
and $\delta_i$.
\end{proof}

\subsection{Remark on nonfree actions}\label{rem:nonfree}
Theorem~\ref{thm:cnilp} can be reformulated as follows:

There exists a constant $\eps(m)>0$ such that if $N^m$  is a Riemannian manifold which admits a free discrete isometric action by a  group $\Gamma$
such that $\sec(N)\z> -\eps(m)$ and $\diam(N/\Gamma)<1$,
then $\Gamma$ is $C(m)$-nilpotent.

As was pointed out to us by B.~Wilking,
in the above reformulation of Theorem~\ref{thm:cnilp}
one can actually remove the assumption
that the $\Gamma$ action be free.

\begin{cor}\label{cor:nonfree}
There exists a constant $\eps(m)>0$ such that
if $N^m$  is a Riemannian manifold
which admits a discrete  isometric action by a  group $\Gamma$
such that $\sec(N)> -\eps(m)$ and $\diam(N/\Gamma)<1$,
then $\Gamma$ is $C(m)$-nilpotent.
\end{cor}
\begin{proof}
Let $\eps=\eps(m)$ be as provided by Theorem~\ref{thm:cnilp}
and suppose $N$ satisfies the assumptions of the corollary for this $\eps$.
Let $F$ be the frame bundle of $N$.
Then the action of $\Gamma$ on $N$ lifts to a free isometric action on $F$.
As was observed in~\cite{FY},
using Cheeger's rescaling trick,
$F$ can be equipped with a $\Gamma$ invariant metric
satisfying $\sec(F)>-\eps(m)$ and $\diam(F/\Gamma)<1$.
Since the induced action of $\Gamma$ on $F$ is free,
the claim of the corollary now follows from Theorem~\ref{thm:cnilp}.
\end{proof}

\section{Proof of the Fibration Theorem}\label{sec:fib}
\subsection{} Let $M$ be an almost nonnegatively curved manifold.
Denote by $M_n\z=(M,g_n)$
a sequence of Riemannian metrics on $M$ such that $\sec(M_n)\ge -1/n$
and $\diam(M_n)\le 1/n$.

Let us denote by $\tilde M$ the universal cover of $M$ and by
$\tilde M_n\to M_n$
the universal Riemannian covering of $M_n$
(that is, $\tilde M$ equipped with the pull back of the metric $g_n$ on $M$).

By \cite{FY}, passing to a finite cover
we may assume that $\Gamma=\pi_1M$ is a nilpotent group without torsion.
Hence, to prove the topological part of Theorem~\ref{intro:main3},
it is sufficient to show the following:

\begin{thm}\label{thm:fib}

Let $M$ be a closed almost nonnegatively curved $m$-manifold
such that $\Gamma=\pi_1M$ is a nilpotent group without torsion.
Then  $M$
is the total space of a fiber bundle
$$
 F\to M\to N
$$
where the base $N$ is a nilmanifold and the fiber $F$ is simply connected.
\end{thm}

The assumption on $\Gamma$ implies that
we can fix a series
$$\Gamma=\Gamma_0\rhd\Gamma_1\rhd\Gamma_2\rhd\ldots\rhd\Gamma_\ell=\{1\}$$
such that $\Gamma_i$ is normal in $\Gamma$ and $\Gamma_i/\Gamma_{i+1}\cong \mathbb{Z}$.

Let us first us give an informal proof.

\subsection{An informal proof of Theorem~\ref{thm:fib}}
We use induction to construct the bundles
$F_i\longrightarrow M\overset{f_i}{\longrightarrow} N_i$, where each $N_i$ is a
nilmanifold with $\pi_1N_i=\Gamma/\Gamma_i$
and $\pi_1F_i\cong \Gamma_i$.
Since the base of induction is trivial,
we are only interested in the induction step.

Fix $p\in N_i$,
and let $F_i(p)$ be the fiber over $p$.
For any sufficiently large $n$
choose a subgroup $G_i=G_i(n)$
such that $\Gamma_{i}\lhd G_i\lhd \Gamma_{i+1}$ and
$[\Gamma_i: G_i]$ is finite,
but sufficiently large so that the cover
$\bar{F}_i(p)$ of $F_i(p)$
corresponding to $G_i$ is Hausdorff close to a unit circle $\mathbb{S}^{1}$.

Construct now a bundle map $\phi_p\colon \bar F_i(p)\to \mathbb{S}^{1}$
by lifting distance functions from $\mathbb{S}^{1}$
(This can be done by a slight generalization
of a construction in ~\cite{FY}  and~\cite{BGP}).
Let $\omega_p=d\phi_p$.

Then $\omega_p$
is a closed integral non-degenerate one-form on $F_i(p)$.
Since deck transformations are isometries,
after  averaging by $\mathbb{Z}_a$, where $a= [\Gamma_i: G_i]$,
we can assume that
$\omega_p$ is $\mathbb{Z}_a$-invariant.
Thus $\omega_p$ descends to a form on $F_i(p)$
which when integrated gives a bundle map $F_i(p)$ onto a small $\mathbb{S}^{1}$.

Note that altho this bundle is defined only up to rotations of $\mathbb{S}^{1}$,
its fibers are well-defined.

Since $\Gamma_{i+1}$ is normal in $\Gamma$,
the choice of the covering $\bar{F}_i(p)$ of $F_i(p)$
is unambiguous for all $p\in N_i$.
By using a partition of unity on $N_i$ we can
glue the forms $\omega_p$ into a global 1-form on $M$
which  satisfies the following properties:
\begin{enumerate}[a)]
\item $\omega|_{F(p)}$ is closed  and integral for any $p$;
\item  $\omega|_{F(p)}$ is non-degenerate.
\end{enumerate}

Integrating $\omega$ over the various $F(p)$'s
we construct a continuous family of bundles $F_p\to \mathbb{S}^{1}$.
The level sets
partition each $F(p)$ and hence the whole $M$ into fibers of a fiber bundle,
whose quotient space is then a circle bundle $N_{i+1}$ over $N_i$.\qed

\medskip

This gives a good idea of the proof.
However, to make it precise some extra work has to be done.
In particular,  one has to be careful with the construction of $\omega$.
To make this construction possible
we have to keep track of how $F(p)$  was obtained.
Namely, we have to use that the fiber $F(p)$
was obtained by a construction as in Yamaguchi's
fibration theorem (see~\cite{Yam} or~\cite{BGP}).
This makes the induction proof quite technical.

We now proceed with the rigorous proof of Theorem~\ref{thm:fib}.

\subsection{Proof of Theorem~\ref{thm:fib}}
Let us denote by  $\tilde M_{n,i}$
the Riemannian covering of $M_n$ with respect to $\Gamma_i$.

For any choice of marked points $p_n$
we have that
$$(\tilde M_{n,i},p_n,\pi_1M)\GHto (\mathbb{R}^i,0,\mathbb{R}^i)$$
in equivariant Gromov--Hausdorff convergence,
where $\mathbb{R}^i$ acts on itself by translations.
Indeed, the limit space must be a
nonnegatively curved simply connected Alexandrov space,
and since $\diam M_n\to0$ we have that it possesses
a transitive group action by a nilpotent group.
Then Euclidean space, acting as a group of translations,
is here the only choice,
and it is easy to see that the dimension of the limit  must be equal to $i$.

Therefore $(\tilde M_{n},p_n,\pi_1M)\GHto (\mathbb{R}^\ell,0,\mathbb{R}^\ell)$,
and we may also assume that for each $i$
we have that $(\tilde M_{n},p_n,\Gamma_i)\GHto (\mathbb{R}^\ell,0, \mathbb{R}^{\ell-i})$,
where $\mathbb{R}^{\ell-i}$ is the coordinate subspace of $\mathbb{R}^\ell$
which corresponds to the first $\ell-i$ elements of the standard basis.

Now, let us give a technical definition:

If $R$ is a Riemannian manifold,
let us denote by $\widetilde\dist_p$ the average of a distance function
over a small ball around $p$.
This enables us
to work with the $C^{1}$  function $\widetilde\dist_p$
instead of the Lipschitz function $\dist_p$.

\begin{defn}\label{def:close}
Let $R_n\GHto R$ be a sequence of Riemannian $m$-manifolds
with curvature $\ge \kappa$ which Gromov--Hausdorff converges
to a Riemannian $m'$-manifold $R$, where $m'\le m$.
A sequence of forms $\omega_n$ on $R_n$ is said to
$\eps$-approximate a form $\omega$ on $R$,
if
\begin{enumerate}[(i)]
\item
for any point $p\in R$ there is a neighborhood $U\ni p$
which admits a distance chart $f\colon U\to \mathbb{R}^{m'}$,
$$f(x)=(\dist_{a_1}(x),\dist_{a_2}(x),\dots,\dist_{a_{m'}}(x))$$
which is a smooth regular map,
and

\item
smooth lifts of $f$ to $U_n\subset R_n$
give, for all large $n$, regular maps
$$f_n(x)=(\widetilde\dist_{a_{1,n}}(x),\widetilde\dist_{a_{2,n}}(x),\dots,
\widetilde\dist_{a_{m',n}}(x))$$
with $a_{i,n}\in M_n$, $a_{i,n}\to a_{n}$
such that
$$|(f_n\circ f^{-1})^*(\omega)-\omega_n|_{C^{0}}<\eps$$
for all sufficiently large $n$.

\end{enumerate}

\end{defn}

Theorem~\ref{thm:fib} now easily follows from the following lemma:

\begin{lem}
Given $\eps>0$ there is a sequence of one-forms
$$\{\omega_{1,n},\omega_{2,n},\cdots,\omega_{k,n}\}$$
on $\tilde M_n$
with the following properties:
\begin{enumerate}[(i)]
\item For each $i$, $\omega_{i,n}$ is a $\pi_1M$-invariant form on $\tilde M_n$.
\item The forms $\omega_{i,n}$ $\eps$-approximate
the coordinate forms $dx_i$ on $\mathbb{R}^k$.
In particular, the forms $\{\omega_{i,n}\}$ are nowhere zero
and almost orthonormal at each point.
\item
If for any $j<i$ it holds that
$\omega_{j,n}(X)=\omega_{j,n}(Y)=0$, then $d\omega_{i,n}(X,Y)=0$.
In particular,
for each $i$ and all sufficiently large $n$,
the distribution corresponding to the system of equations
$$\omega_{j,n}(X)=0 \quad\text{for all}\quad j\le i$$
defines  on $\tilde M_n$ a foliation $\mathcal F_{i,n}$.

\item
If $\tilde F_{i,n}(x)\subset \tilde M_n$ denotes
the fiber of the foliation $\mathcal F_{i,n}$ containing the point $x\in \tilde M_n$,
then each $\tilde F_{i,n}(x)$ is $\Gamma_i$-invariant;
that is,  for any $\gamma\in \Gamma_i$
one has that $\tilde F_{i,n}(x)=\tilde F_{i,n}(\gamma x)$.
 Moreover, the action of $\Gamma_i$ on $\tilde F_{i,n}(x)$ is cocompact
for each $i$. In particular, $\mathcal F_{i,n}$ induces on $M_n$
the structure of a fiber bundle.
\end{enumerate}
\end{lem}

\begin{proof}
We will construct these forms by induction.
Assume that we have already constructed one-forms
$\omega_1,\omega_2,\ldots, \omega_{i-1}$
which meet all the required properties.
They give a $\pi_1M$-invariant fibration of $\tilde M_n$ by submanifolds
$\tilde F_{i-1,n}(x)$ thru each point $x\in \tilde M_n$, with
tangent spaces defined by the equations $\omega_j(X)\z=0$ for $j=1,\ldots,i-1$.

Denote by $\theta\colon \mathbb{R}\to [0,1]$
a smooth monotone function
which is equal to $1$ before $0$ and $0$ after $1$.
Choose numbers $\delta_n>0$ slowly converging to $0$,
and let $\Theta_{i,n}\colon \tilde M_n\z\to \mathbb{R}_+$ be the function defined by
$$\Theta_{i,n}(x)
=\min_{y\in F_{i-1,n}(x)}\{\theta(|p_n y|/\delta_n)\}.$$
Clearly $\Theta_{i,n}$ is a continuous $\Gamma_{i-1}$-invariant function
which is constant on each $F_{i-1,n}(x)$.
Moreover, for large $n$,
$\Theta_{i,n}$ has support in  some $C_i{\cdot} \delta_n$-neighborhood
of $F_{i-1,n}(p_n)$, and is equal to $1$
in some $c_i{\cdot} \delta_n$-neighborhood of $F_{i-1,n}(p_n)$.

Now let $\phi\colon \mathbb{R}\to [0,1]$
be a smooth nondecreasing function
which is $0$ before $1/2$ and $1$ after $3/2$.
Consider the form
$$\omega_{i,n}'=\Theta_{i,n}{\cdot} d(\phi\circ \widetilde\dist_{\Gamma_{i}a_{i,n}}),$$
where $a_{i,n}\in \tilde M_n$ is a sequence of points
converging to $ -e_i\in \mathbb{R}^\ell$,
and $\widetilde\dist_{\Gamma_{i}a_{i,n}}$ is the
average of $\dist_{\Gamma_{i}x}$ for $x$ in a small ball around $a_{i,n}$.
The support of $\omega_i'$ has two components,
one which  contains $ p_n$ (notice here that $p_n\to 0\in \mathbb{R}^\ell$),
and another which does not.
(It follows from the construction
that  the limit of $F_{i-1,n}(p_n)$ is a coordinate plane  in $\mathbb{R}^\ell$).

Set $\omega_{i,n}'':=\omega_{i,n}'$ on the component of $ p_n$,
and let this form be $0$ otherwise.
Clearly, $\omega_{i,n}''$ is then a continuous $\Gamma_{i}$-invariant form
whose restriction to ${\tilde F_{i-1,n}(x)}$ is exact.
Moreover, each level set of its integral over
${\tilde F_{i-1,n}(x)}$ is $\Gamma_{i}$-invariant.

By construction, the form $\omega''/|\omega''|$ is now
 (in the sense of definition~\ref{def:close})
close to $dx_i$  at the points where
$|\omega''|\not=0$.
Take
$$\omega_{i,n}=c\sum_{\gamma\in \Gamma/\Gamma_i} \gamma\omega',$$
where the coefficient $c$ is chosen in such a way that  $|\omega_{i,n}(p_n)|=1$.
As $\delta_n$ is a sequence slowly converging to zero,
we may assume that $\diam (M_n)/\delta_n\to 0$.
Therefore, $\omega_{i,n}$ is the form we need.
\end{proof}

Notice that the proof actually shows
that the fibers in Theorem~\ref{thm:fib} are almost nonnegatively curved manifolds
in the generalized sense with $k=\ell$. Therefore,
 the proof of Theorem~\ref{intro:main3} is complete. \hfill$\qed$

\section{Open questions}
We would like to conclude this work by posing a number of related open questions.

\subsection{Is the torsion contained in the center?}
As was noted earlier, Theorem ~\ref{intro:main2} is new even for manifolds of nonnegative curvature.
For such manifolds it is known that their fundamental groups  are almost abelian,
and Fukaya and Yamaguchi conjectured the following  (see \cite{FY}):

\begin{conj}[Fukaya--Yamaguchi]\label{con:c-ab}
The fundamental group of a nonnegatively curved  $m$-manifold is $C(m)$-abelian.
\end{conj}

In this regard we would like to pose the following two conjectures:

\begin{mconj}\label{con:tor}
There is $C=C(m)$ such that if $M^m$ is almost nonnegatively curved then there is a nilpotent subgroup $N\subset \pi_1M$ of index $\le C$ whose torsion is contained in its center (or, at least, whose torsion is commutative).
\end{mconj}

\begin{conj}\label{pi_2}
If $M^m$ is almost nonnegatively curved,
then the action of $\pi_1M$ on $\pi_2M$ is almost trivial, or maybe even $C(m)$-trivial;
that is, there exists a finite index subgroup of $\pi_1M$
(or, respectively,  a subgroup of index $\le C(m)$)
which acts trivially on $\pi_2M$.
\end{conj}


Conjecture ~\ref{con:tor} implies in particular that the fundamental
groups of closed positively curved $m$-manifolds are $C(m)$-abelian.

In fact, as was pointed out to us by B.~Wilking, if true,
Conjecture~\ref{con:tor} would also imply a positive answer to Conjecture~\ref{con:c-ab}.
Indeed, if $\sec(M)\ge0$, then the universal cover $\tilde M$ of $M$ is isometric to the product $\mathbb{R}^n\times K$, where $K$ is a compact Riemannian manifold and the $\pi_1M$ action  on $\mathbb{R}^n\times K$ is diagonal.
It follows from \cite[Cor. 6.3]{wilking} that one can deform the metric on $M$ so that its universal cover is still isometric to $\mathbb{R}^n\times K$ and the induced action on $K$ is finite.
By passing, as in the proof of Corollary~\ref{cor:nonfree}, to the induced action on the frame bundle of $K$, one reduces the statement to Conjecture~\ref{con:tor}.

We tried to prove these conjectures by studying successive blow-ups of the
collapsing sequence $M_n$ as done in Section \ref{sec:constr}.

We can prove  Conjectures~\ref{con:tor} and \ref{pi_2} in the case where  all spaces $A_i$ which appear in the construction in  Section \ref{sec:constr} are closed Riemannian manifolds; 
see \cite{KPT}.
Moreover, we believe we have an argument to prove it if all $A_i$'s are Alexandrov spaces without boundary.

It seems that if we would have just a slightly better understanding of collapsing to a ray, then we could prove the conjectures.
Here is the simplest related question which we cannot solve:

\medskip

\begin{quest}
Let $M_n=(\mathbb{S}^2\times\mathbb{R}^2,g_n)$ be a sequence of complete Riemannian manifolds with  $\sec(M_n)\ge -\eps_n$, where $\eps_n\to0$ as $n\to\infty$.
Assume that for a sequence of points $p_n\in M_n$ we have that $(M_n,p_n)\GHto (\mathbb{R}_+,0)$.
Let $q_n\in M_n$  be a sequence of points such that $|p_n q_n|=1$ and such that there is a sequence of rescalings $\lambda_n\to\infty$ such that
$(\lambda_n{\cdot} M_n,q_n)\GHto \mathbb{S}^2\times \mathbb{S}^{1}\times \mathbb{R}$, where the latter space is equipped with the product of the canonical metrics.
\begin{enumerate}[(i)]

\item Can it happen that $(\lambda_n{\cdot} M_n,p_n)\GHto(\mathbb{R}_+,0)$?

\item Is it true that  the dimension of the Gromov--Hausdorff limit of $(\lambda_n {\cdot} M_n,p_n)$ is at least 3?

\item What are the possible limits of $(\lambda_n{\cdot}  M_n,p_n)$?

\end{enumerate}

\end{quest}
\medskip

Conjectures \ref{con:c-ab} and  \ref{con:tor} are also related to the following conjecture of Rong
(cf. \cite{Ro2,Ro1}):

\begin{conj}[Rong] \label{con:rong} Positively curved $m$-manifolds have
$C(m)$-cyclic fundamental groups.
\end{conj}

This conjecture has been proved by Rong \cite{Ro2} under the
additional assumption of a uniform upper curvature bound. We also
believe that if one could carry out the above program for proving
Conjecture~\ref{con:tor}, one would have a good shot at handling
Rong's Conjecture as well.

\subsection{The simply connected case}
So far we have only discussed manifolds with nontrivial fundamental groups.
However, some of our arguments  also work in a more general setting.
We hope that it might be possible  to use them to obtain new restrictions on simply connected almost nonnegatively curved manifolds as well as on collapsing with a lower curvature bound.

Let us indicate one possible approach to do so.

Let us denote by $\cM(F)$ the space of self homotopy equivalences of a manifold $F$.
Assume now that $F$ is simply connected and that $\tilde f\colon \mathbb{S}^{k}\times F\to F$ is a map such that $\tilde f_u\colon F\to F$ is homotopic to the identity for some (and therefore ANY) $u\in \mathbb{S}^{k}$.
Then $\tilde f$ represents an element $\alpha=[\tilde f]\in \pi_{k}(\cM(F))$.
Let $g$ be a Riemannian metric on $F$.
Define
$$\dil_g(\tilde f)=\max_{u\in \mathbb{S}^k}\dil_g(\tilde f_u),$$
where  $\dil_g(\tilde f_u)$ stands for the optimal Lipschitz constant of $\tilde f_u$ with respect to $g$.
For any $\alpha\in \pi_k(\cM(F))$ define
$$\dil_g(\alpha)=\inf_{[h]=\alpha}\dil_g(h).$$
Finally, define
$$\DIL(\alpha)=\inf_g\dil_g(\alpha)$$
over all Riemannian metrics $g$ on $F$ and
$$\DIL_+(\alpha)=\inf_g\dil_g(\alpha)$$
over all Riemannian metrics $g$ on $F$ with $\diam(F,g)\le 1$ and $\sec(g)\ge -1$.

Clearly, both $\DIL(\alpha)$ and $\DIL_+(\alpha)$ are  homotopy invariants of $\alpha$.

Now suppose that $M_n\GHto \mathbb{S}^{k+1}$ is a sequence of Riemannian manifolds collapsing to a round sphere with $\sec(M_n)\ge k$.
By Yamaguchi's fibration theorem, we have that $M_n$ is a fiber bundle over $\mathbb{S}^{k+1}$ with almost nonnegatively curved fiber $F_n$.
This bundle is classified by an element $\alpha$  of $\pi_k(Aut(F_n))$ and by using our gradient flow technique we can estimate $\DIL_+(\alpha)$ (and hence $\DIL(\alpha)$) from above.

Therefore, if one could find examples of a simply connected $F$ and an $\alpha$ with arbitrary big $\DIL_+(\alpha)$, one would obtain new restrictions on collapsing to a sphere with curvature bounded from below, and probably more restrictions for the topological type of manifolds with lower curvature and upper diameter bounds in general.
In fact, $F$ need not be simply connected as long as the total space of the bundle $F\to M\to \mathbb{S}^{k+1}$ is.

While we believe that finding examples with arbitrary large  $\DIL(\alpha)$) is very difficult
(and might even be impossible), we have several candidates to produce large $\DIL_+(\alpha)$.

On the other hand, the problem of finding examples of $\alpha$ with
$\DIL(\alpha)>1$
seems quite interesting in its own right and might have other applications unrelated to collapsing.

\medskip

Let us next describe  some possible sources of  examples with
$\DIL_+(\alpha)>1$:

\begin{ex}
Obviously, if dil$_g(h)=1$, then $h_u$ is a homotopy of isometries of $(F,g)$.
Let $G$ be the isometry group of $F$.
Then $G$ can be viewed as a subset of $\cM(F)$.
Therefore, if $[h]\not= 0$ in $\pi_k(\cM(F))$, then $[h_u]\not=0$ in $\pi_kG$.
Now $G$ is a compact Lie group, in particular, $\pi_2G=0$
(and even more generally $\pi_{2{\cdot} n}G$ is finite).
On the other hand, there are spaces $F$ such that the space $\cM(F)$ might have nontrivial second homotopy;
for example, the canonical metric on $F\z=\mathrm{SU}(6)/(\mathrm{SU}(3)\z\times \mathrm{SU}(3))$ has nonnegative curvature,
and it follows
from  \cite[Chapter 5, Example 4.14]{OT},
that
$\pi_2(\cM(F))\otimes \mathbb{Q}$ is nontrivial.
Therefore, there is  an $\alpha \in \pi_2(\cM(F))$ such that
dil$_g(\alpha)>1$ for any metric $g $ on $F$;
we believe it should be true that $\DIL_+(\alpha)>1$.
Still, it might happen that $\DIL(\alpha)=1$.
\end{ex}

Another possible source of such manifolds is provided by the following example due to D.~Sullivan.
\begin{ex}
Let $N^7$ be the total space of an $\mathbb{S}^3$ bundle over $\mathbb{S}^4$ with zero Euler class and nontrivial $p_1$.
Clearly $N^7$ is rationally equivalent to $\mathbb{S}^4\times \mathbb{S}^3$.
In particular, its minimal model has no nontrivial derivations of degree $-1$.
Therefore,  by~\cite[13.3]{Su}, there exists a diffeomorphism $f\colon N\to N$ which is homotopic to the identity but such that the obstruction to it being diffeotopic to the identity is a nonzero element of $H^3(N,\mathbb{Z})\cong \mathbb{Z}$.
Let $M^8$ be the mapping cylinder of $f$.
Clearly $M$ is homotopy equivalent to $N\times \mathbb{S}^{1}$ and hence it is spin with signature zero.
On the other hand, by construction, $p_1^2(M)\ne 0$.
Since the signature of $M$ is zero we must necessarily have that $p_2(M)\ne 0$  and hence $\hat{A}(M)\ne 0$.
In particular, by the Atiyah--Hirzebruch theorem, $M$ does not admit an $\mathbb{S}^{1}$ action  and hence the corresponding element $\alpha \in \pi_1(\cM(M))$ has dil$_g(\alpha)>1$ for any metric $g$ on $M$.
\end{ex}

\begin{rmk}
As was mentioned in the introduction, it is known that a spin
manifold $X$ of almost nonnegative Ricci curvature has
$\hat{A}(X)\le 2^{\dim X/2}$ (\cite[page 41]{G5}, \cite{Ga}).
Clearly, a finite cover of the manifold $M$ constructed above
violates this restriction and therefore $M$ does not admit almost
nonnegative Ricci curvature. However,  it could possibly be almost
nonnegatively curved in the generalized sense.
\end{rmk}

\subsection{Further questions}

Recall that a simply connected space $C$ is called  \emph{rationally
elliptic} if it is homotopy equivalent to a finite CW-complex and

\[
\dim [\pi_*(C,\mathbb{Q})]<\infty.
\]

A conjecture of Grove--Halperin~\cite{GH} says that simply connected nonnegatively curved manifolds are rationally elliptic. 
This conjecture was extended by Grove to include almost nonnegatively curved manifolds~\cite{Grv}.
Later, Totaro has proposed the following definition of  rationally elliptic spaces which covers manifolds with infinite fundamental groups:

A connected topological space $X$ is \emph{rationally elliptic}
if it is homotopy equivalent to a finite CW complex,
it has a finite covering which is a nilpotent space
and its universal cover is rationally elliptic in the ordinary sense.

With this definition one can extend Grove's conjecture  
to non simply connected manifolds as follows:

\begin{conj}\label{conj-grove}
Any almost nonnegatively curved manifold in the generalized sense is rationally elliptic.
\end{conj}

Theorem ~\ref{intro:main1}  reduces this conjecture to the simply connected case which is undoubtedly the most difficult part of the problem.

It has been shown in \cite{PaPe} that if $M$ is a nilpotent closed manifold which admits a Riemannian metric with zero topological entropy, then its universal cover $\tilde{M}$ is rationally elliptic.
Coupled with  Theorem ~\ref{intro:main1} this means that to prove Conjecture~\ref{conj-grove} it would be sufficient to show that a manifold with almost nonnegative curvature in
the
generalized sense admits a metric with zero topological entropy.
However, we think that this
 might be wrong in general.

\medskip

As was pointed out in the discussion in the Introduction before Theorem~\ref{intro:main3},  it already follows from Yamaguchi's fibration theorem
and \cite{FY}
that a finite cover of an almost nonnegatively curved manifold maps onto a nilmanifold
  with homotopy fiber a simply connected closed manifold.
While this is formally weaker than the statement of Theorem~\ref{intro:main3}, it would be interesting to
have an answer to the following, purely topological, question:

\begin{quest}
Let $ M \overset{f}{\longrightarrow} N$ be a map from a  closed manifold $M$
to  a nilmanifold $N$ such that  the homotopy fiber of $f$ is  a simply connected closed manifold.
Is it true that after passing to a finite cover, the map
$f$ becomes homotopic to a fiber bundle projection?
\end{quest}

\begin{quest}
Is it true  that  manifolds which are almost nonnegatively curved in the generalized sense
are almost nonnegatively curved?
\end{quest}

In view of Theorems ~\ref{intro:main1} and ~\ref{intro:main2} it is
also reasonable to pose the following question:

\begin{quest}
Is it true that almost nonnegatively curved $m$-manifolds $M^m$
are $C(m)$-nilpotent spaces?
\end{quest}

It is clear from the proof of Theorems ~\ref{intro:main1} and ~\ref{intro:main2}
that this is true if the universal cover
of $M^m$ has torsion free integral cohomology.

In view of Theorem~\ref{intro:main2} it is moreover natural to raise
the following question:
\begin{quest}
Can one give an explicit bound on $C(m)$ in Theorem~\ref{intro:main2}?
\end{quest}

\small
\bibliographystyle{alpha}


\end{document}